\documentclass{article}%
\usepackage{graphicx}
\usepackage[intlimits]{amsmath}
\usepackage{latexsym}
\usepackage{amsfonts}
\usepackage{amssymb}%
\newtheorem{theorem}{Theorem}

\newtheorem{corollary}[theorem]{Corollary}

\newtheorem{lemma}[theorem]{Lemma}

\newenvironment{proof}[1][Proof]{\noindent{\textbf {#1}  }}  {\hfill$\Box$\bigskip}

\begin{document}

\title{Joints in graphs}
\author{B\'{e}la Bollob\'{a}s\thanks{Department of Mathematical Sciences, University
of Memphis, Memphis TN 38152, USA} \thanks{Trinity College, Cambridge CB2 1TQ,
UK} \thanks{Research supported in part by NSF grant Itr\# DARPA grant
F33615-01-C-1900.} \ and Vladimir Nikiforov$^{\ast}$}
\maketitle

\begin{abstract}
In 1969 Erd\H{o}s proved that if $r\geq2$ and $n>n_{0}\left(  r\right)  ,$
every graph $G$ of order $n$ and $e\left(  G\right)  >t_{r}\left(  n\right)  $
has an edge that is contained in at least $n^{r-1}/\left(  10r\right)  ^{6r}$
$\left(  r+1\right)  $-cliques. In this note we improve this bound to
$n^{r-1}/r^{r+5}.$ We also prove a corresponding stability result.


\end{abstract}

\section{Introduction}

Our notation and terminology are standard (see, e.g. \cite{Bol98}). Thus,
$G\left(  n\right)  $ is a graph of order $n$ and $G\left(  n,m\right)  $ is a
graph of order $n$ and size $m;$ for a graph $G$ and a vertex $u\in V\left(
G\right)  $ we write $\Gamma\left(  u\right)  $ for the neighborhood of $u$;
$d_{G}\left(  u\right)  =\left\vert \Gamma\left(  u\right)  \right\vert $ is
the degree of $u;$ we write $d\left(  u\right)  $ for $d_{G}\left(  u\right)
$ when there is no danger of confusion. We denote by $k_{r}\left(  G\right)  $
the number of $r$-cliques of $G.$ We let $T_{r}\left(  n\right)  $ be the
Tur\'{a}n graph of order $n$ with $r$ classes and set $t_{r}\left(  n\right)
=e\left(  T_{r}\left(  n\right)  \right)  $.

Erd\H{o}s \cite{Erd69} proved that if $r\geq2$ and $n>n_{0}\left(  r\right)
,$ every graph $G=G\left(  n,t_{r}\left(  n\right)  +1\right)  $ contains at
least $n^{r-1}/\left(  10r\right)  ^{6r}$ cliques of order $\left(
r+1\right)  $ sharing an edge. He used this result to estimate the minimum
number of cliques in certain graphs.

In this note we strengthen and extend this result of Erd\H{o}s. We start with
a general definition. Let $p,q,r$ be integers with $p\geq r,$ $q>r\geq1.$ We
call the union of a $p$-clique $H$ and $t$ $q$-cliques, each one intersecting
$H$ in exactly $r$ vertices, \emph{a }$\left(  p,q,r\right)  $\emph{-joint }of
size $t$\emph{ }and denote it by $J_{t}^{\left(  p,q,r\right)  }.$ The maximum
size of a $\left(  p,q,r\right)  $-joint in a graph $G$ is called the $\left(
p,q,r\right)  $\emph{-jointsize }of\emph{ }$G$ and is denoted by $js^{\left(
p,q,r\right)  }\left(  G\right)  .$

Observe that, in general, there may be many nonisomorphic $\left(
p,q,r\right)  $-joints with the same parameters $p,q,$ $r.$

In terms of joints the above assertion of Erd\H{o}s can be stated as follows:
for every integer $r\geq2$ and $n>n_{0}\left(  r\right)  ,$
\begin{equation}
js^{\left(  2,r+1,2\right)  }\left(  G\left(  n,t_{r}\left(  n\right)
+1\right)  \right)  \geq\frac{n^{r-1}}{\left(  10r\right)  ^{6r}}.
\label{Erdb}%
\end{equation}

In this note we shall show that, in fact, if $r\geq2$ and $n>r^{8}$ then
\[
js^{\left(  2,r+1,2\right)  }\left(  G\left(  n,t_{r}\left(  n\right)
+1\right)  \right)  >\frac{n^{r-1}}{r^{r+5}}.
\]

Moreover, we shall show that if $r\geq2,$ $n>r^{8},$ and $0<\alpha<36r^{-8}$
then, for every graph $G=G\left(  n\right)  $ with $e\left(  G\right)
>t_{r}\left(  n\right)  -\alpha n^{2},$ either
\[
js^{\left(  2,r+1,2\right)  }\left(  G\left(  n,t_{r}\left(  n\right)
+1\right)  \right)  >\left(  1-\frac{1}{r^{3}}\right)  \frac{n^{r-1}}{r^{r+5}}%
\]
or $G$ contains an induced $r$-chromatic subgraph of order at least $\left(
1-2\sqrt{\alpha}\right)  n.$

\section{Preliminary results}

Recall that the following basic properties of the Tur\'{a}n graph
$T_{r}\left(  n\right)  $
\begin{equation}
\delta\left(  T_{r}\left(  n\right)  \right)  =\left\lfloor \frac{r-1}%
{r}n\right\rfloor \label{turdel}%
\end{equation}
and%
\begin{equation}
t_{r}\left(  n\right)  =t_{r}\left(  n-1\right)  +\delta\left(  T_{r}\left(
n\right)  \right)  . \label{turec}%
\end{equation}
Furthermore,%
\begin{equation}
t_{r}\left(  n\right)  =\frac{r-1}{2r}\left(  n^{2}-t^{2}\right)  +\binom
{t}{2}, \label{turval}%
\end{equation}
where $t$ is the remainder of $n$ modulo $r,$ and so%
\begin{equation}
\frac{r-1}{2r}n^{2}-\frac{r}{8}\leq t_{r}\left(  n\right)  \leq\frac{r-1}%
{2r}n^{2}. \label{turest}%
\end{equation}

\subsection{Bounds on $k_{r+1}\left(  G\right)  $ and $js^{\left(
2,r+1,2\right)  }\left(  G\right)  $}

We start by establishing lower bounds for $k_{r+1}\left(  G\right)  $ and
$js^{\left(  2,r+1,2\right)  }\left(  G\right)  $ in a graph $G$ of order $n$
with $e\left(  G\right)  >\frac{r-1}{2r}n^{2}.$

Note that $N$ cliques $K_{r+1}$ of a graph $G$ cover some edge at least
$N\binom{r+1}{2}/e\left(  G\right)  $ times, and so%
\begin{equation}
js^{\left(  2,r+1,2\right)  }\left(  G\right)  >k_{r+1}\left(  G\right)
\binom{r+1}{2}\binom{n}{2}^{-1}. \label{prel}%
\end{equation}

\begin{lemma}
\label{leNSMM}For all $r\geq3,$ $c>0$, if $G=G\left(  n\right)  $ and%
\begin{equation}
e\left(  G\right)  >\left(  \frac{r-1}{2r}+c\right)  n^{2} \label{cond1}%
\end{equation}
then%
\begin{equation}
k_{r+1}\left(  G\right)  >2c\frac{r}{r+1}\left(  \frac{n}{r}\right)  ^{r+1}
\label{lok}%
\end{equation}
and
\begin{equation}
js^{\left(  2,r+1,2\right)  }\left(  G\right)  >2c\left(  \frac{n}{r}\right)
^{r-1} \label{loj}%
\end{equation}

\end{lemma}

\begin{proof}
In \cite{MoMo62} Moon and Moser stated the following assertion whose complete
proof apparently appeared for the first time in \cite{KhNi78} (see also
\cite{Lov79}, Problem 11.8).

If $G=G\left(  n\right)  $ and $k_{s}\left(  G\right)  >0$ then
\[
\frac{\left(  s+1\right)  k_{s+1}\left(  G\right)  }{sk_{s}\left(  G\right)
}-\frac{n}{s}\geq\frac{sk_{s}\left(  G\right)  }{\left(  s-1\right)
k_{s-1}\left(  G\right)  }-\frac{n}{s-1}.
\]
Equivalently, if $q$ is the clique number of $G$ then, for $q>s>t\geq1,$ we
have%
\begin{equation}
\frac{\left(  s+1\right)  k_{s+1}\left(  G\right)  }{sk_{s}\left(  G\right)
}-\frac{n}{s}\geq\frac{\left(  t+1\right)  k_{t+1}\left(  G\right)  }%
{tk_{t}\left(  G\right)  }-\frac{n}{t}. \label{MoMo}%
\end{equation}

Since Tur\'{a}n's theorem and (\ref{cond1}) imply $k_{r+1}\left(  G\right)
>0,$ setting $t=1$ in (\ref{MoMo}), we find that
\[
\frac{\left(  s+1\right)  k_{s+1}\left(  G\right)  }{sk_{s}\left(  G\right)
}-\frac{n}{s}\geq\frac{2e\left(  G\right)  }{n}-n>\left(  -\frac{1}%
{r}+2c\right)  n
\]
for every $s=2,...,r.$ Hence,
\[
\frac{\left(  s+1\right)  k_{s+1}\left(  G\right)  }{sk_{s}\left(  G\right)
}>\left(  \frac{1}{s}-\frac{1}{r}+2c\right)  n
\]
for every $s=1,...,r.$ Multiplying these inequalities for $s=1,...,r,$ we find
that%
\[
\frac{\left(  r+1\right)  k_{r+1}\left(  G\right)  }{n}\geq n^{r}\prod
_{s=1}^{r}\left(  \frac{1}{s}-\frac{1}{r}+2c\right)  >2cn^{r}\prod_{s=1}%
^{r-1}\left(  \frac{1}{s}-\frac{1}{r}\right)  =\frac{2c}{r^{r}}n^{r},
\]
and hence (\ref{lok}) holds.

Taking into account (\ref{prel}), we find that%
\[
js^{\left(  2,r+1,2\right)  }\left(  G\right)  \geq\binom{r+1}{2}%
k_{r+1}\left(  G\right)  \binom{n}{2}^{-1}>\binom{r+1}{2}\frac{2c}{\left(
r+1\right)  r^{r}}n^{r+1}\binom{n}{2}^{-1},
\]
and (\ref{loj}) follows.
\end{proof}

Since the inequality $2e\left(  G\right)  \geq\delta\left(  G\right)  v\left(
G\right)  $ holds for every graph $G,$ Lemma \ref{leNSMM} implies the
following corollary.

\begin{corollary}
\label{cor} For all $r\geq3,$ $c>0$, if $G=G\left(  n\right)  $ and%
\[
\delta\left(  G\right)  >\left(  \frac{r-1}{r}+c\right)  n
\]
then%
\[
k_{r+1}\left(  G\right)  >c\frac{r}{r+1}\left(  \frac{n}{r}\right)  ^{r+1}%
\]
and
\[
js^{\left(  2,r+1,2\right)  }\left(  G\right)  >c\left(  \frac{n}{r}\right)
^{r-2}.
\]

\end{corollary}

\subsection{A Bonferroni-Zarankievicz type inequality}

Suppose $r\geq3,$ $X$ is a set of cardinality $n,$ and $A_{1},...,A_{r}$ are
subsets of $X.$ For every $k\in\left[  r\right]  ,$ set%
\[
S_{k}=\sum_{1\leq i_{1}<...<i_{k}\leq r}\mu\left(  A_{i_{1}}\cap...\cap
A_{i_{k}}\right)  .
\]
Then the following theorem holds.

\begin{theorem}
\label{TYPMS} If $1\leq k\leq r$ then%
\begin{equation}
S_{k}\geq\binom{\left\lfloor S_{1}/n\right\rfloor }{k-1}\left(  S_{1}%
-\frac{k-1}{k}\left(  \left\lfloor S_{1}/n\right\rfloor +1\right)  n\right)
\label{goal}%
\end{equation}

\end{theorem}

\begin{proof}
Let $H$ be a bipartite graph whose color classes are the sets $\left[
r\right]  $ and $X$, and $i\in\left[  r\right]  $ is joined to $u\in X$ iff
$u\in A_{i}.$ Clearly,%
\begin{equation}
S_{1}=e\left(  H\right)  =\sum_{u\in X}d_{H}\left(  u\right)  \label{scon}%
\end{equation}
\ and%
\[
S_{k}=\sum_{u\in X}\binom{d_{H}\left(  u\right)  }{k}.
\]
The convexity of $\binom{x}{k}$ implies that the minimum of $S_{k},$ subject
to (\ref{scon}), is attained when every vertex $u$ has degree $d_{H}\left(
u\right)  =\left\lfloor S_{1}/n\right\rfloor $ or $d_{H}\left(  u\right)
=\left\lceil S_{1}/n\right\rceil .$ Letting $l$ be the number of those $u$
with $d_{H}\left(  u\right)  =\left\lceil S_{1}/n\right\rceil $ and setting
$x=l/n,$ we see that%
\[
\left(  1-x\right)  \left\lfloor S_{1}/n\right\rfloor +x\left\lceil
S_{1}/n\right\rceil =S_{1}/n,
\]
and so, $x=S_{1}/n-\left\lfloor S_{1}/n\right\rfloor .$ Since%
\[
\binom{\left\lfloor S_{1}/n\right\rfloor }{k-1}=\binom{\left\lfloor
S_{1}/n\right\rfloor +1}{k}-\binom{\left\lfloor S_{1}/n\right\rfloor }{k}%
\]
for $x>0,$ we have%
\begin{align*}
S_{k}  &  \geq\left(  1-x\right)  n\binom{\left\lfloor S_{1}/n\right\rfloor
}{k}+xn\binom{\left\lceil S_{1}/n\right\rceil }{k}\\
&  =n\binom{\left\lfloor S_{1}/n\right\rfloor }{k}+xn\left(  \binom
{\left\lfloor S_{1}/n\right\rfloor +1}{k}-\binom{\left\lfloor S_{1}%
/n\right\rfloor }{k}\right) \\
&  =n\binom{\left\lfloor S_{1}/n\right\rfloor }{k}+n\left(  S_{1}%
/n-\left\lfloor S_{1}/n\right\rfloor \right)  \binom{\left\lfloor
S_{1}/n\right\rfloor }{k-1}\\
&  =\binom{\left\lfloor S_{1}/n\right\rfloor }{k-1}\left(  S_{1}-\frac{k-1}%
{k}\left(  \left\lfloor S_{1}/n\right\rfloor +1\right)  n\right)  ,
\end{align*}
as claimed.
\end{proof}

\begin{lemma}
\label{lebonf1} Suppose $r\geq2,$ $0<a<1/r\left(  r+1\right)  ,$ $X$ is a set
of cardinality $n,$ and $A_{1},...,A_{r+1}$ are subsets of $X.$ If%
\[
\sum_{i=1}^{r+1}\left\vert A_{i}\right\vert \geq\left(  r-\frac{1}{r}-\left(
r+1\right)  a\right)  n.
\]
Then some two members of $\left\{  A_{1},...,A_{r+1}\right\}  $ have at least%
\[
\left(  \frac{r-2}{r}+\frac{2}{r^{2}\left(  r+1\right)  }-\frac{2\left(
r-1\right)  }{r}a\right)  n
\]
elements in common.
\end{lemma}

\begin{proof}
Applying Theorem \ref{TYPMS} with $k=2$ to the sets $A_{1},...,A_{r+1},$ we
find that
\begin{align*}
S_{2}  &  \geq\binom{r-1}{1}\left(  S_{1}-\frac{r}{2}n\right)  \geq\left(
r-1\right)  \left(  r-\frac{1}{r}-\left(  r+1\right)  a-\frac{r}{2}\right)
n\\
&  =\left(  \frac{r\left(  r-1\right)  }{2}-\frac{r-1}{r}-\left(
r^{2}-1\right)  a\right)  n.
\end{align*}
Since there are $\binom{r+1}{2}$ pairwise intersections $A_{i}\cap A_{j},$ for
some $1\leq k<l\leq r+1$ we have%
\begin{align*}
\left\vert A_{k}\cap A_{l}\right\vert  &  \geq S_{2}\binom{r+1}{2}^{-1}%
\geq\left(  \frac{r\left(  r-1\right)  }{2}-\frac{r-1}{r}-\left(
r^{2}-1\right)  a\right)  \binom{r+1}{2}^{-1}n\\
&  =\left(  \frac{r-1}{r+1}-\frac{2\left(  r-1\right)  }{r^{2}\left(
r+1\right)  }-\frac{2\left(  r-1\right)  }{r}a\right)  n\\
&  =\left(  \frac{r-2}{r}+\frac{2}{r^{2}\left(  r+1\right)  }-\frac{2\left(
r-1\right)  }{r}a\right)  n.
\end{align*}

\end{proof}

The idea of the following lemma is due to Erd\H{o}s; our proof techniques
allow to improve his bound considerably.

\begin{lemma}
\label{leKd} Suppose $r\geq3.$ If a graph $G=G\left(  n\right)  $ contains a
$K_{r+1}$ and%
\[
\delta\left(  G\right)  >\left(  \frac{r-1}{r}-\frac{1}{r^{2}\left(
r^{2}-1\right)  }\right)  n
\]
then
\[
js^{\left(  2,r+1,2\right)  }\left(  G\right)  >\frac{n^{r-1}}{r^{r+3}}.
\]

\end{lemma}

\begin{proof}
Indeed, let $U$ be the vertex set of an $\left(  r+1\right)  $-clique in $G.$
Then%
\[
\sum_{i\in U}\left\vert A_{i}\right\vert >\left(  r+1\right)  \left(
\frac{r-1}{r}-\frac{1}{r^{2}\left(  r^{2}-1\right)  }\right)  n=\left(
r-\frac{1}{r}-\frac{r+1}{r^{2}\left(  r^{2}-1\right)  }\right)  n.
\]
Hence, by Lemma \ref{lebonf1}, there are distinct $u,v\in U$ such that
$M=\left\vert \Gamma\left(  u\right)  \cap\Gamma\left(  v\right)  \right\vert
$ satisfies
\begin{align}
\left\vert M\right\vert  &  \geq\left(  \frac{r-2}{r}+\frac{2}{r^{2}\left(
r+1\right)  }-\frac{2\left(  r-1\right)  }{r}\frac{1}{r^{2}\left(
r^{2}-1\right)  }\right)  n\nonumber\\
&  =\left(  \frac{r-2}{r}+\frac{2\left(  r-1\right)  }{r^{3}\left(
r+1\right)  }\right)  n. \label{lom}%
\end{align}

For the graph $G\left[  M\right]  $ induced by the set $M$ we have%
\begin{align}
\delta\left(  G\left[  M\right]  \right)   &  \geq\delta\left(  G\right)
-\left(  n-\left\vert M\right\vert \right)  >\left(  \frac{r-1}{r}-\frac
{1}{r^{2}\left(  r^{2}-1\right)  }\right)  n-\left(  n-\left\vert M\right\vert
\right) \nonumber\\
&  =\left\vert M\right\vert -\left(  \frac{1}{r}+\frac{1}{r^{2}\left(
r^{2}-1\right)  }\right)  n. \label{lodm}%
\end{align}

By routine calculations we find that, for $r\geq3,$%
\[
\left(  \frac{1}{r-2}-\frac{1}{r^{2}\left(  r-1\right)  ^{2}}\right)  \left(
\frac{r-2}{r}+\frac{2\left(  r-1\right)  }{r^{3}\left(  r+1\right)  }\right)
>\frac{1}{r}+\frac{1}{r^{2}\left(  r^{2}-1\right)  }.
\]

Recalling (\ref{lom}), this implies%
\[
\left\vert M\right\vert \left(  \frac{1}{r-2}-\frac{1}{r^{2}\left(
r-1\right)  ^{2}}\right)  >\left(  \frac{1}{r}+\frac{1}{r^{2}\left(
r^{2}-1\right)  }\right)  n,
\]
and furthermore,%
\[
\left\vert M\right\vert -\left(  \frac{1}{r}+\frac{1}{r^{2}\left(
r^{2}-1\right)  }\right)  n>\left(  \frac{r-3}{r-2}+\frac{1}{r^{2}\left(
r^{2}-1\right)  }\right)  \left\vert M\right\vert .
\]
Hence, from inequality (\ref{lodm}) we see that
\[
\delta\left(  G\left[  M\right]  \right)  \geq\left(  \frac{r-3}{r-2}+\frac
{1}{r^{2}\left(  r-1\right)  ^{2}}\right)  \left\vert M\right\vert .
\]
In view of (\ref{lom}), Corollary \ref{cor} implies%
\[
k_{r-1}\left(  G\left[  M\right]  \right)  \geq\frac{r-2}{r^{2}\left(
r-1\right)  ^{3}}\left(  \frac{\left\vert M\right\vert }{r-2}\right)
^{r-1}>\frac{r-2}{r^{2}\left(  r-1\right)  ^{3}}\frac{n^{r-1}}{r^{r-1}}%
>\frac{n^{r-1}}{r^{r+3}}.
\]

To complete the proof observe that the number of $\left(  r+1\right)
$-cliques of $G$ containing the edge $uv$ is exactly $k_{r-1}\left(  G\left[
M\right]  \right)  .$
\end{proof}

\section{\label{Jexst}Existence of large joints $J^{\left(  2,r+1,2\right)  }%
$}

In this section we shall prove a Tur\'{a}n type result for large joints as
stated in Theorem \ref{Turj} below. We start with the following technical result.

\begin{theorem}
\label{Thexj} If $r\geq2$ and $n>r^{8}$, every graph $G=G\left(  n\right)  $
with
\begin{equation}
e\left(  G\right)  >t_{r}\left(  n\right)  \label{cnd}%
\end{equation}
has an induced subgraph $G^{\prime}=G\left(  n^{\prime}\right)  $ with
$n^{\prime}>\left(  1-1/r^{2}\right)  n$\ such that either%
\begin{equation}
K_{r+1}\subset G^{\prime}\text{, \ \ \ and \ \ \ }\delta\left(  G^{\prime
}\right)  >\left(  \frac{r-1}{r}-\frac{1}{r^{2}\left(  r^{2}-1\right)
}\right)  n^{\prime}, \label{prop3}%
\end{equation}
or
\begin{equation}
e\left(  G^{\prime}\right)  >\left(  \frac{r-1}{2r}+\frac{1}{r^{4}\left(
r^{2}-1\right)  }\right)  \left(  n^{\prime}\right)  ^{2}. \label{prop4}%
\end{equation}

\end{theorem}

\begin{proof}
Let the sequence $u_{1},...,u_{n}$ be an enumeration of the vertices such that
$d\left(  u_{1}\right)  =\delta\left(  G\right)  \ $and
\[
d\left(  u_{i}\right)  =\delta\left(  G-u_{1}-...-u_{i-1}\right)  \text{ \ for
}1<i\leq n.
\]

Set $G_{0}=G,$ and set $G_{i}=G-u_{1}-...-u_{i},$ $i=1,...,n-1,$ so that
\begin{equation}
e\left(  G_{i}\right)  -e\left(  G_{i+1}\right)  =\delta\left(  G_{i}\right)
\label{eqn}%
\end{equation}
for every $i\in\left[  n-1\right]  .$

Set $\beta=\frac{1}{r^{2}\left(  r^{2}-1\right)  }$ and let $k-1$ be the
largest integer such that $1\leq k\leq n$ and
\[
\delta\left(  G_{k-1}\right)  \leq\left(  \frac{r-1}{r}-\beta\right)  \left(
n-k+1\right)  .
\]

From (\ref{eqn}), for every $s\in\left[  k\right]  ,$ we have%
\begin{align*}
e\left(  G\right)  -e\left(  G_{s}\right)   &  =\sum_{i=0}^{s-1}\delta\left(
G_{i}\right)  \leq\left(  \frac{r-1}{r}-\beta\right)  \sum_{i=0}^{s-1}\left(
n-i\right)  \\
&  \leq\left(  \frac{r-1}{r}-\beta\right)  \left(  \binom{n+1}{2}%
-\binom{n-s+1}{2}\right)  \\
&  <\left(  \frac{r-1}{r}-\beta\right)  \left(  \frac{n^{2}}{2}-\frac{\left(
n-s\right)  ^{2}}{2}+\frac{s}{2}\right)  .
\end{align*}
From (\ref{cnd}) and (\ref{turest}) we have%
\[
e\left(  G\right)  >\frac{r-1}{2r}n^{2}-\frac{r}{8}.
\]
Hence, for every $s\in\left[  k\right]  $, we deduce%
\begin{align}
e\left(  G_{s}\right)   &  >e\left(  G\right)  -\left(  \frac{r-1}{r}%
-\beta\right)  \left(  \frac{n^{2}}{2}-\frac{\left(  n-s\right)  ^{2}}%
{2}+\frac{s}{2}\right)  \nonumber\\
&  >\frac{r-1}{2r}n^{2}-\frac{r}{8}-\left(  \frac{r-1}{2r}-\frac{\beta}%
{2}\right)  \left(  n^{2}-\left(  n-s\right)  ^{2}+s\right)  \nonumber\\
&  =\beta\frac{n^{2}}{2}+\left(  \frac{r-1}{2r}-\frac{\beta}{2}\right)
\left(  n-s\right)  ^{2}-\left(  \frac{r-1}{2r}-\frac{\beta}{2}\right)
s-\frac{r}{8}\nonumber\\
&  >\frac{r-1}{2r}\left(  n-s\right)  ^{2}+\frac{\beta}{2}\left(
n^{2}-\left(  n-s\right)  ^{2}\right)  -\frac{r}{8}-\frac{s}{2}.\label{los}%
\end{align}

In the rest of the proof we shall consider two cases - \emph{(a) }$k>n/r^{2}$
and \emph{(b) }$k\leq n/r^{2}.$

\emph{(a) }Let\emph{ }$n>r^{8},$ assume that $k>n/r^{2},$ and set
$l=\left\lfloor n/r^{2}\right\rfloor .$ Then we have%
\begin{equation}
n-l\leq\left(  1-\frac{1}{r^{2}}\right)  \left(  n+1\right)  , \label{ord}%
\end{equation}
implying%
\[
\left(  n-l\right)  ^{2}\leq\left(  n+1\right)  ^{2}\left(  1-\frac{1}{r^{2}%
}\right)  ^{2}\leq\left(  n^{2}-\frac{n}{\beta}\right)  \left(  1+\frac
{2}{r^{2}}\right)  ^{-1}.
\]
Hence, from (\ref{los}), it follows
\begin{align*}
e\left(  G_{l}\right)   &  >\frac{r-1}{2r}\left(  n-l\right)  ^{2}+\frac
{\beta}{2}\left(  n^{2}-\left(  n-l\right)  ^{2}\right)  -\frac{r}{8}-\frac
{n}{2r^{2}}\\
&  >\frac{r-1}{2r}\left(  n-l\right)  ^{2}+\frac{\beta}{2}\left(  n^{2}%
-\frac{n}{\beta}-\left(  n-l\right)  ^{2}\right) \\
&  >\left(  \frac{r-1}{2r}+\frac{\beta}{r^{2}}\right)  \left(  n-l\right)
^{2}.
\end{align*}
This, together with (\ref{ord}), implies (\ref{prop4}) with $G^{\prime}%
=G_{l}.$

\emph{(b) }Assume that $k\leq n/r^{2}$. The way the graphs $G_{1},...,G_{k}$
are constructed, together with (\ref{turec}) and (\ref{cnd}), implies%
\[
e\left(  G_{k}\right)  >t_{r}\left(  n-k\right)  ,
\]
and by Tur\'{a}n's theorem $K_{r+1}\subset G.$ Since $\left(  n-k\right)
\geq\left(  1-1/r^{2}\right)  n$, and
\[
\delta\left(  G_{k}\right)  >\left(  \frac{r-1}{r}-\frac{1}{r^{2}\left(
r^{2}-1\right)  }\right)  \left(  n-k\right)  ,
\]
condition (\ref{prop3}) holds with $G^{\prime}=G_{k}.$ The proof is completed.
\end{proof}

After this proposition we are ready to prove our main theorem, strengthening
inequality (\ref{Erdb}).

\begin{theorem}
\label{Turj} For $r\geq2$ and $n>r^{8}$, every graph $G=G\left(  n\right)  $
with $e\left(  G\right)  \geq t_{r}\left(  n\right)  $ satisfies%
\begin{equation}
js^{\left(  2,r+1,2\right)  }\left(  G\right)  >\frac{n^{r-1}}{r^{r+5}}
\label{ourb}%
\end{equation}
unless $G=T_{r}\left(  n\right)  .$
\end{theorem}

\begin{proof}
Assume first that $e\left(  G\right)  >t_{r}\left(  n\right)  .$ By Theorem
\ref{Thexj} $G$ contains an induced subgraph $G^{\prime}=G\left(  n^{\prime
}\right)  $ with $n^{\prime}>\left(  1-1/r^{2}\right)  n$\ and such that
either (\ref{prop3}) or (\ref{prop4}) holds. If (\ref{prop3}) is true,
applying Lemma \ref{leKd} to the graph $G^{\prime},$ we see that%
\[
js^{\left(  2,r+1,2\right)  }\left(  G^{\prime}\right)  \geq\frac{\left(
n^{\prime}\right)  ^{r-1}}{r^{r+3}}>\left(  1-\frac{1}{r^{2}}\right)
^{r-1}\frac{n^{r-1}}{r^{r+3}}>\left(  1-\frac{1}{r}\right)  \frac{n^{r-1}%
}{r^{r+3}},
\]
and the assertion follows.

If (\ref{prop4}) holds then, by Lemma \ref{leNSMM}, we see that%
\begin{align*}
js^{\left(  2,r+1,2\right)  }\left(  G^{\prime}\right)   &  \geq\frac{2}%
{r^{4}\left(  r^{2}-1\right)  }\left(  \frac{n^{\prime}}{r}\right)
^{r-1}>\frac{2}{r^{4}\left(  r^{2}-1\right)  }\left(  1-\frac{1}{r^{2}%
}\right)  ^{r-1}\left(  \frac{n}{r}\right)  ^{r-1}\\
&  >\frac{2}{r^{4}\left(  r^{2}-1\right)  }\left(  1-\frac{1}{r}\right)
\left(  \frac{n}{r}\right)  ^{r-1}>\frac{n^{r-1}}{r^{r+5}},
\end{align*}
and the assertion follows.

Assume now that $e\left(  G\right)  =t_{r}\left(  n\right)  .$ If $G$ has a
vertex $u$ with $d\left(  u\right)  <\delta\left(  T_{r}\left(  n\right)
\right)  $ then
\[
e\left(  G-u\right)  >t_{r}\left(  n-1\right)  ,
\]
and therefore, the graph $G-u$ contains an induced subgraph $G^{\prime
}=G\left(  n^{\prime}\right)  $ with $n^{\prime}>\left(  1-1/r^{2}\right)
\left(  n-1\right)  $\ and such that either (\ref{prop3}) or (\ref{prop4})
holds. Using the arguments from the first part of our proof we see that
either
\begin{align*}
js^{\left(  2,r+1,2\right)  }\left(  G^{\prime}\right)   &  >\frac{\left(
n^{\prime}\right)  ^{r-1}}{r^{r+5}}>\left(  1-\frac{1}{r^{2}}\right)
^{r-1}\left(  1-\frac{1}{r^{8}}\right)  ^{r-1}\frac{n^{r-1}}{r^{r+3}}\\
&  >\left(  1-\frac{r-1}{r^{2}}\right)  \left(  1-\frac{r-1}{r^{8}}\right)
\frac{n^{r-1}}{r^{r+3}}>\frac{n^{r-1}}{r^{r+5}},
\end{align*}
or%
\begin{align*}
js^{\left(  2,r+1,2\right)  }\left(  G^{\prime}\right)   &  \geq\frac{2}%
{r^{4}\left(  r^{2}-1\right)  }\left(  \frac{n^{\prime}}{r}\right)
^{r-1}>\frac{2}{r^{4}\left(  r^{2}-1\right)  }\left(  1-\frac{1}{r^{2}%
}\right)  ^{r-1}\left(  \frac{n-1}{r}\right)  ^{r-1}\\
&  >\frac{2}{r^{4}\left(  r^{2}-1\right)  }\left(  1-\frac{1}{r}\right)
\left(  1-\frac{1}{r^{8}}\right)  ^{r-1}\left(  \frac{n}{r}\right)
^{r-1}>\frac{n^{r-1}}{r^{r+5}},
\end{align*}
completing the proof in this case.

It remains the case when $\delta\left(  G\right)  =\delta\left(  T_{r}\left(
n\right)  \right)  .$ Hence, in view of $n>r^{8},$ we find that
\begin{equation}
\delta\left(  G\right)  =\left\lfloor \frac{r-1}{r}n\right\rfloor \geq\left(
\frac{r-1}{r}-\frac{1}{r^{2}\left(  r^{2}-1\right)  }\right)  n. \label{deg1}%
\end{equation}
If $G\neq T_{r}\left(  n\right)  ,$ Tur\'{a}n's theorem implies that $G$
contains a $K_{r+1};$ thus, in view of Lemma \ref{leKd} and (\ref{deg1}), the
proof is completed.
\end{proof}

Note that (\ref{ourb}) is tight up to a factor of order at most $r^{-6},$ as
seen by taking the graph $T_{r}\left(  n\right)  $ and adding an edge to its
largest chromatic class.

\section{A stability theorem about large joints $J^{\left(  2,r+1,2\right)  }%
$}

Theorem \ref{Turj} may be used to prove a stability result about large joints
$J^{\left(  2,r+1,2\right)  }$ as stated in the theorem below. In the course
of our proof we shall need the following result of Andr\'{a}sfai, Erd\H{o}s
and S\'{o}s \cite{AES74}: if $G$ is a $K_{r+1}$-free graph of order $n$ with
minimal degree
\[
\delta\left(  G\right)  >\left(  1-\frac{3}{3r-1}\right)  n
\]
then $G$ is $r$-chromatic.

\begin{theorem}
\label{stabj} Let $r\geq2,$ $n>r^{8}$, and $0<\alpha<r^{-8}/36.$ If a graph
$G=G\left(  n\right)  $ satisfies
\[
e\left(  G\right)  >\left(  \frac{r-1}{2r}-\alpha\right)  n^{2},
\]
then either
\begin{equation}
js^{\left(  2,r+1,2\right)  }\left(  G\right)  >\left(  1-\frac{1}{r^{3}%
}\right)  \frac{n^{r-1}}{r^{r+5}}, \label{minjs}%
\end{equation}
or $G$ contains an induced $r$-chromatic subgraph $G_{0}$ of order at least
$\left(  1-2\sqrt{\alpha}\right)  n$ with minimum degree
\begin{equation}
\delta\left(  G_{0}\right)  >\left(  1-\frac{1}{r}-6\sqrt{\alpha}\right)  n.
\label{mindg}%
\end{equation}

\end{theorem}

\begin{proof}
We may assume that $\alpha n^{2}\geq1,$ since otherwise we have $e\left(
G\right)  \geq t_{r}\left(  n\right)  $ and the assertion follows from Theorem
\ref{Turj}. Set
\begin{equation}
\varepsilon=2\sqrt{\alpha}<\frac{1}{3r^{4}}, \label{maxeps}%
\end{equation}
and define $M_{\varepsilon}\subset V$ as%
\[
M_{\varepsilon}=\left\{  u\in V\left(  G\right)  :d\left(  u\right)
\leq\left(  \frac{r-1}{r}-\varepsilon\right)  n\right\}  .
\]

Assume that (\ref{minjs}) does not hold. Our aim is to show that \emph{(a)
}$\left\vert M_{\varepsilon}\right\vert <2\varepsilon n,$ and \emph{(b)} the
subgraph $G_{0}$ of $G$ induced by $V\left(  G\right)  \backslash
M_{\varepsilon}$ has the properties required in the theorem.%

\begin{align*}
\frac{\left(  s+1\right)  k_{s+1}\left(  G\right)  }{sk_{s}\left(  G\right)
}-\frac{n}{s}  &  >\frac{3k_{3}\left(  G\right)  }{2m}-\frac{n}{2}\geq\frac
{1}{2m}\sum_{u\in V\left(  G\right)  }d^{2}\left(  u\right)  -n\\
&  =\frac{2m}{n}+\frac{1}{2m}\sum_{u\in V\left(  G\right)  }\left(  d\left(
u\right)  -\frac{2m}{n}\right)  ^{2}-n\\
&  \geq\frac{2m}{n}-n+\frac{1}{2m}\sum_{u\in M_{\varepsilon}}\varepsilon
^{2}>\frac{2m}{n}-n+\varepsilon^{2}\left\vert M_{\varepsilon}\right\vert
\end{align*}

\emph{(a) }Assume, for a contradiction, that $\left\vert M_{\varepsilon
}\right\vert \geq2\varepsilon n$ and let $M^{\prime}\subset M_{\varepsilon}$
satisfy
\begin{equation}
\left(  1-\sqrt{1/2}\right)  \varepsilon n<\left\vert M^{\prime}\right\vert
<\left(  1+\sqrt{1/2}\right)  \varepsilon n. \label{bnds}%
\end{equation}
Such a set $M^{\prime}$ exists since $\sqrt{2}\varepsilon n>2\sqrt{2}%
\sqrt{\alpha}n>2\sqrt{2}.$ Let $G^{\prime}$ be the subgraph of $G$ induced by
$V\backslash M^{\prime}.$ Then%
\begin{align*}
e\left(  G\right)   &  =e\left(  G^{\prime}\right)  +e\left(  M^{\prime
},V\backslash M^{\prime}\right)  +e\left(  M^{\prime}\right)  \leq e\left(
G^{\prime}\right)  +\sum_{u\in M^{\prime}}d\left(  u\right) \\
&  \leq e\left(  G^{\prime}\right)  +\left\vert M^{\prime}\right\vert \left(
\frac{r-1}{r}-\varepsilon\right)  n.
\end{align*}

Observe that the second inequality of (\ref{bnds}) implies that%
\[
n-\left\vert M^{\prime}\right\vert >n-2\varepsilon n.
\]
Hence, if
\[
e\left(  G^{\prime}\right)  >\frac{r-1}{2r}\left(  n-\left\vert M^{\prime
}\right\vert \right)  ^{2}%
\]
then, by Theorem \ref{Turj} and (\ref{maxeps}),
\begin{align*}
js^{\left(  2,r+1,2\right)  }\left(  G\right)   &  \geq js^{\left(
2,r+1,2\right)  }\left(  G^{\prime}\right)  >\frac{\left(  n-\left\vert
M^{\prime}\right\vert \right)  ^{r-1}}{r^{r+5}}>\left(  1-2\varepsilon\right)
^{r-1}\frac{n^{r-1}}{r^{r+5}}\\
&  >\left(  1-2\left(  r-1\right)  \varepsilon\right)  ^{r}\frac{n^{r-1}%
}{r^{r+5}}>\left(  1-\frac{1}{r^{3}}\right)  \frac{n^{r-1}}{r^{r+5}}.
\end{align*}
Thus (\ref{minjs}) holds, contradicting our assumption.

Consequently we may assume that%
\[
e\left(  G^{\prime}\right)  \leq\frac{r-1}{2r}\left(  n-\left\vert M^{\prime
}\right\vert \right)  ^{2}.
\]
Since%
\begin{align*}
e\left(  G^{\prime}\right)   &  \geq e\left(  G\right)  -\sum_{u\in M}d\left(
u\right) \\
&  \geq\left(  \frac{r-1}{2r}-\alpha\right)  n^{2}-\left\vert M^{\prime
}\right\vert \left(  \frac{r-1}{r}-\varepsilon\right)  n,
\end{align*}
it follows that
\[
\frac{r-1}{2r}\left(  n-\left\vert M^{\prime}\right\vert \right)  ^{2}%
\geq\left(  \frac{r-1}{2r}-\alpha\right)  n^{2}-\left\vert M^{\prime
}\right\vert \left(  \frac{r-1}{r}-\varepsilon\right)  n.
\]

Setting $x=\left\vert M^{\prime}\right\vert /n$ we find that
\[
\frac{r-1}{2r}\left(  1-x\right)  ^{2}+x\left(  \frac{r-1}{r}-\varepsilon
\right)  -\left(  \frac{r-1}{2r}-\alpha\right)  \geq0.
\]
and so,%
\[
x^{2}-2\varepsilon x+2\alpha\geq0.
\]
Hence, either%
\[
\left\vert M^{\prime}\right\vert \leq\left(  \varepsilon-\sqrt{\varepsilon
^{2}-2\alpha}\right)  n=\varepsilon\left(  1-\sqrt{1/2}\right)  n
\]
or%
\[
\left\vert M^{\prime}\right\vert \geq\left(  \varepsilon+\sqrt{\varepsilon
^{2}-2\alpha}\right)  n=\varepsilon\left(  1+\sqrt{1/2}\right)  n,
\]
contradicting (\ref{bnds}). Therefore, $\left\vert M_{\varepsilon}\right\vert
<2\varepsilon n.$

\emph{(b) }Note first that $G_{0}$ has $n-\left\vert M_{\varepsilon
}\right\vert >\left(  1-2\sqrt{\alpha}\right)  n$ vertices. By our choice of
$M_{\varepsilon},$ for $u\in V\backslash M_{\varepsilon},$ we have
\begin{equation}
d_{G}\left(  u\right)  >\left(  \frac{r-1}{r}-\varepsilon\right)  n,
\label{mindg1}%
\end{equation}
so
\begin{equation}
d_{G_{0}}\left(  u\right)  >\left(  \frac{r-1}{r}-\varepsilon\right)
n-\left\vert M_{\varepsilon}\right\vert >\left(  \frac{r-1}{r}-3\varepsilon
\right)  n=\left(  \frac{r-1}{r}-6\sqrt{\alpha}\right)  n, \label{mindg2}%
\end{equation}
and (\ref{mindg}) holds.

All that remains to prove is that $G_{0}$ is $r$-chromatic. From
(\ref{mindg2}) we have%
\begin{align}
\delta\left(  G_{0}\right)   &  >\left(  \frac{r-1}{r}-6\sqrt{\alpha}\right)
n\geq\left(  \frac{r-1}{r}-6\sqrt{\alpha}\right)  v\left(  G_{0}\right)
\nonumber\\
&  >\left(  \frac{r-1}{r}-\frac{1}{r^{4}}\right)  v\left(  G_{0}\right)
>\left(  1-\frac{3}{3r-1}\right)  v\left(  G_{0}\right)  \label{mindeg2}%
\end{align}
If $G_{0}$ contains a $K_{r+1}$, by Lemma \ref{leKd} we have%
\begin{align*}
js^{\left(  2,r+1,2\right)  }\left(  G\right)   &  \geq js^{\left(
2,r+1,2\right)  }\left(  G_{0}\right)  >\frac{\left(  n-\left\vert M^{\prime
}\right\vert \right)  ^{r-1}}{r^{r+5}}>\left(  1-2\varepsilon\right)
^{r-1}\frac{n^{r-1}}{r^{r+5}}\\
&  >\left(  1-2\left(  r-1\right)  \varepsilon\right)  \frac{n^{r-1}}{r^{r+5}%
}>\left(  1-\frac{1}{r^{3}}\right)  \frac{n^{r-1}}{r^{r+5}}.
\end{align*}
Therefore, (\ref{minjs}) holds, contradicting our assumption.

We may assume that $G_{0}$ is $K_{r+1}$-free. In view of (\ref{mindeg2}), the
theorem of Andr\'{a}sfai, Erd\H{o}s and S\'{o}s implies that $G_{0}$ is
$r$-chromatic, completing our proof.
\end{proof}

\end{document}